\documentclass[MSNbibl,number,citesort,dvips]{arxbj}
\usepackage{upgreek}
\usepackage{graphicx}
\aid{0}
\volume{20}
\issue{2}
\pubyear{2014}
\firstpage{775}
\lastpage{802}
\doi{10.3150/13-BEJ506} 

\makeatletter
\newcommand{\RMe}{\mathrm{e}}

\newcommand{\mrmd}{\,\mathrm{d}}

\newtheorem{lemma}{Lemma}[section]
\newremark{example}{Example}[section]
\newtheorem{theorem}{Theorem}[section]
\newtheorem{corollary}{Corollary}[section]
\newremark{remark}{Remark}[section]

\newcommand{\R}{\mathbb R}
\newcommand{\N}{\mathbb N}

\makeatother

\begin{document}
\begin{frontmatter}

\title{Maximum likelihood characterization of distributions}
\runtitle{Maximum likelihood characterization}

\begin{aug}
\author[1]{\fnms{Mitia} \snm{Duerinckx}\thanksref{1,e1}\ead[label=e1,mark]{mitia.duerinckx@ulb.ac.be}},
\author[1]{\fnms{Christophe} \snm{Ley}\thanksref{1,e2}\ead[label=e2,mark]{christophe.ley@ulb.ac.be}} \and
\author[2]{\fnms{Yvik} \snm{Swan}\corref{}\thanksref{2}\ead[label=e3]{yvik.swan@uni.lu}}
\runauthor{M. Duerinckx, C. Ley and Y. Swan} 
\address[1]{D\'epartement de Math\'ematique, Universit\'e libre de Bruxelles,
Campus Plaine CP 210, Boulevard du triomphe, B-1050 Bruxelles,
Belgique.\\ \printead{e1,e2}}
\address[2]{Universit\'e du Luxembourg, Campus Kirchberg, UR en math\'ematiques,
6, rue Richard Couden\-hove-Kalergi, L-1359 Luxembourg,
Grand-Duch\'e de Luxembourg. \printead{e3}}
\end{aug}

\received{\smonth{9} \syear{2012}}
\revised{\smonth{12} \syear{2012}}

%
\begin{abstract}
A famous characterization theorem due to C.F. Gauss states that the
maximum likelihood estimator (MLE) of the parameter in a location
family is the sample mean for all samples of all sample sizes if and
only if the family is Gaussian. There exist many extensions of this
result in diverse directions, most of them focussing on location and
scale families. In this paper, we propose a unified treatment of
this literature by providing general MLE characterization theorems
for one-parameter group families (with particular attention on
location and scale parameters). In doing so, we provide tools for
determining whether or not a given such family is
MLE-characterizable, and, in case it is, we define the fundamental
concept of \textit{minimal necessary sample size} at which a given
characterization holds. Many of the cornerstone references on this
topic are retrieved and discussed in the light of our findings, and
several new characterization theorems are provided. Of particular
interest is that one part of our work, namely the introduction of
so-called \textit{equivalence classes} for MLE characterizations, is a
modernized version of Daniel Bernoulli's viewpoint on maximum
likelihood estimation.
%
\end{abstract}

%

\end{frontmatter}

\section{Introduction}\label{intro}

In probability and statistics, a characterization theorem occurs
whenever a given law or a given class of laws is the only one which
satisfies a certain property. While probabilistic characterization
theorems are concerned with distributional aspects of (functions~of)
random variables, statistical characterization theorems rather deal
with properties of \textit{statistics}, that is, measurable functions of
a set of independent random variables \mbox{(observations)} following a
certain distribution. Examples of probabilistic characterization
theorems include:
\begin{enumerate}[-]
\item[-] Stein-type characterizations of probability laws, inspired from
the classical Stein and Chen characterizations of the normal and the
Poisson distributions, see Stein \cite{S72},
Chen \cite{C75}; 

\item[-] maximum entropy characterizations, see Cover and Thomas \cite
{cover2006elements},
Chapter 11;
\item[-] conditioning characterizations of probability laws
such as, for example, determining the marginal distributions of the random
components $X$ and $Y$ of the vector $(X,Y)'$ if only the
conditional distribution of $X$ given $X+Y$ is known, see
Patil and Seshadri~\cite{PS64}.
\end{enumerate}
The class of statistical characterization theorems includes:
\begin{enumerate}[-]
\item[-] characterizations of probability distributions by means of order
statistics, see Galambos \cite{G72} and Kotz \cite{K74} for an overview
on the
vast literature on this subject;
\item[-] Cram\'{e}r-type characterizations, see Cram{\'e}r \cite{Cr36};
\item[-] characterizations of probability laws by means of one linear
statistic, of identically distributed statistics or of the
independence between two statistics, see Lukacs \cite{L56}.
\end{enumerate}
Besides their evident mathematical interest per se,
characterization theorems also provide a better understanding of the
distributions under investigation and sometimes offer unexpected
handles to innovations which might not have been uncovered
otherwise. For instance, Chen and Stein's characterizations are at
the heart of the celebrated \textit{Stein's method}, see the recent
Chen, Goldstein and Shao \cite{CGS11} and Ross \cite{R12} for an
overview; maximum entropy
characterizations are closely related to the development of important
tools in information theory (see Akaike \cite{A77,A78}),
Bayesian probability theory (see Jaynes~\cite{jaynes1957information}) or
even econometrics (see Wu \cite{wu2003calculation}, Park and Bera \cite{PB09};
characterizations of probability distributions by means of
order statistics are, according to Teicher~\cite{T61}, ``harbingers of
[\ldots]
characterization theorems'', and have been
extensively studied around the middle of the twentieth century by the
likes of Kac, Kagan, Linnik, Lukacs or Rao;
Cram\'{e}r-type characterizations are currently the object of active
research, see Bourguin and Tudor \cite{BoTu11}. This list is by no
means exhaustive and
for further information
about characterization theorems as a whole, we refer to
the extensive and still relevant monograph Kagan, Linnik and Rao \cite
{KLR73}, or to
Kotz \cite{K74}, Bondesson \cite{B97}, Haikady \cite{N06} and the
references therein.

In this paper, we focus on a family of characterization theorems which
lie at the
intersection between probabilistic and statistical characterizations,
the so-called 
\textit{MLE characterizations}.

\subsection{A brief history of MLE
characterizations}\label{sechistory}

We call MLE characterization the characterization of a (family of)
probability distribution(s) via the structure of the Maximum Likelihood
Estimator (MLE) of
a certain parameter of interest (location, scale, etc.).

The first occurrence of such a theorem is in Gauss \cite{G09},
where Gauss showed that the normal (a.k.a. the Gaussian) is the only
location family for which the sample mean $\bar{x}=n^{-1}\sum
_{i=1}^nx_i$ is ``always'' the MLE
of the location parameter. More specifically, Gauss \cite{G09} proved
that, in a location family $g(x-\theta)$ with differentiable density
$g$, the
MLE for $\theta$ is the sample mean for all samples $\mathbf
x^{(n)}=(x_1,\ldots,x_n)$
of all sample sizes $n$ if, and only if, $g(x) = \kappa_{\lambda}
\RMe^{-\lambda x^{2}/2}$
for $\lambda>0$ some constant and $\kappa_{\lambda}$ the adequate
normalizing constant.
Discussed as early as in Poincar{{\'e}} \cite{P12}, this important
result, which we
will throughout call \textit{Gauss' MLE characterization}, has
attracted much attention over the past century
and has spawned a spree of papers about MLE characterization
theorems, the main contributions
(extensions and improvements on different levels, see below) being
due to Teicher \cite{T61}, Ghosh and Rao \cite{GR71}, Kagan, Linnik and
Rao \cite{KLR73}, Findeisen \cite{F82}, Marshall and Olkin \cite{MO93}
and Azzalini and Genton \cite{AG07}. (See also H{\"u}rlimann \cite
{Hu98} for an alternative
approach to
this topic.) For more
information on Gauss' original argument, we refer the reader to the
accounts of Hald \cite{Ha98}, pages 354 and 355, and Chatterjee \cite
{C03}, pages
225--227. See Norden \cite{N72} or the fascinating Stigler \cite{S07}
for an interesting discussion on MLEs and the origins of this
fundamental concept.

The successive refinements of Gauss' MLE characterization contain improvements
on two distinct levels. Firstly, several authors have worked
towards weakening the regularity assumptions on the class of
distributions considered; for instance, Gauss requires
differentiability of $g$, while Teicher \cite{T61} only
requires continuity.
Secondly, many authors have aimed at lowering
the sample size necessary for the characterization to hold (i.e., the
``always''-statement); for instance, Gauss
requires that the sample mean be MLE for the location parameter for
\textit{all} sample sizes simultaneously, Teicher \cite{T61} only
requires that
it be
MLE for samples of sizes 2 and 3 at the same time, while Azzalini and
Genton \cite{AG07}
only need that it be so for a single fixed sample size $n\geq
3$. Note that Azzalini and Genton \cite{AG07} also construct explicit
examples of
non-Gaussian distributions for which the sample mean is the MLE of the
location parameter for the sample size $n=2$.
We already draw the reader's attention to the fact that Azzalini and
Genton's \cite{AG07}
result does not supersede Teicher's \cite{T61}, since the former
require more stringent
conditions (namely differentiability of the $g$'s) than the latter. We
will provide a more detailed discussion on this interesting fact again
at the end of Section \ref{secloc}.


Aside from these ``technical'' improvements, the literature
on MLE characterization theorems also contains evolutions in
different directions which have resulted in a new stream of research,
namely that of discovering new (i.e., different from Gauss'
MLE characterization) MLE characterization theorems. These can be
subdivided into
two categories. On the one hand, MLE characterizations with respect to
the location
parameter but for other forms than the sample mean have been shown to
hold for densities other than the Gaussian. On the other hand, MLE
characterizations
with respect to other parameters of interest than the location
parameter have also been
considered. Teicher \cite{T61} shows that if, under some
regularity assumptions and for \textit{all} sample sizes, the MLE for
the scale
parameter of a scale target distribution is the
sample mean $\bar{x}$ then the target is the exponential distribution,
while if it corresponds to the square root of the sample
arithmetic mean of squares $(\frac1n\sum_{i=1}^nx_i^2)^{1/2}$, then
the target is the standard normal distribution.
Following suit on Teicher's
work, Kagan, Linnik and Rao \cite{KLR73} establish that the sample
median is the MLE for the
location parameter for all samples of size $n=4$ if and only if the
parent distribution is the Laplace law. Also, in
Ghosh and Rao \cite{GR71}, it is shown that there exist distributions
other than the
Laplace for which the sample median at $n=2$ or $n=3$ is MLE.
Ferguson \cite{F62} generalizes Teicher's \cite{T61}
location-based MLE characterization from the Gaussian to a
one-parameter generalized normal distribution, and Marshall and Olkin
\cite{MO93}
generalize Teicher's \cite{T61} scale-based MLE characterization of the
exponential distribution to the gamma distribution with shape
parameter $\alpha>0$ by replacing $\bar{x}$ as MLE for the scale parameter
with $\bar{x}/\alpha$.

There also exist contributions by Buczolich and Sz{\'e}kely \cite{BS89}
where they
investigate situations in which a
weighted average of ordered sample elements can be an MLE of the
location parameter. MLE characterizations in the multivariate setup
have been proposed, inter alia, in Marshall and Olkin \cite
{MO93} and Azzalini and Genton \cite{AG07}. In parallel to all these ``linear''
MLE characterizations there have also been a number of developments
regarding MLE characterizations for spherical distributions,
that is distributions taking their values only on the unit
hypersphere in higher dimensions, see Duerinckx and Ley \cite{DL12} and
the references therein.

Finally, there exists a different stream of MLE
characterization research, inspired by Poincar{{\'e}}~\cite{P12}, in
which one
relaxes the assumptions made on the role of the parameter $\theta$ and
seeks to understand the conditions under which the MLE for this
$\theta$ is an arithmetic mean. This approach to MLE characterization
theorems is known as \textit{Gauss' principle}; quoting Campbell \cite{C70}
``\textit{by Gauss' principle we shall mean that a distribution should
be chosen so that the maximum likelihood estimate of the parameter
$\theta$ is the same as the arithmetic mean estimate given by
$n^{-1}\sum_{k=1}^nT(x_k)$}'', with $\theta= {\mathrm E}(T(X))$ for
some known function $T$. When $T$ is the identity function and one
considers one-parameter location families then we recover the MLE
characterization problem solved by Gauss \cite{G09}. For general
one-parameter families the corresponding problem was solved in
Poincar{{\'e}} \cite{P12}. Campbell \cite{C70} broadens
Poincar\'e's conclusion to general functions $T$ and to the multivariate
setup, while 
Bondesson \cite{B97} further extends
the above works to exponential families with nuisance
parameters. As pointed out by their authors, many of these results
remain valid for discrete distributions; for other MLE
characterizations of discrete probability laws, we also refer to Puig
\cite{P03}
and Puig and Valero \cite{PV06}.

\subsection{Applications of MLE characterizations}
\label{secappl-mle-char}

Perhaps the most remarkable application of MLE characterizations is to
be found in the very origins of this field of study, whose
forefathers sought to define families of probability distributions
which were
``natural'' for a given important problem. 
For instance, the Gaussian distribution was
uncovered by Gauss through his effort of finding the location family
for which the sample mean $\bar{x}$ is a most probable value for
$\theta$, the location parameter. Similarly, Poincar\'{e} in his
``Calcul des Probabilit\'{e}s'' (2nd ed., 1912) derived in Chapter 10
(pages 147--168 in the 1896 edition) a particular case of the
\textit{exponential families of distributions} (see Lehmann and Casella
\cite{LC98}, Section
1.5) by asking for which distributions $\bar{x}$ is the MLE of
$\theta$, without specifying the role of the parameter. We refer to
the historical remarks at the end of Bondesson \cite{B97} and the references
therein for more information on both the works of Gauss and
Poincar\'{e}. Following Gauss' ideas, von Mises \cite{vm18} defined the circular
(i.e., spherical in dimension two) analogue of the Gaussian
distribution by looking for the circular distribution whose circular
location parameter always has the circular sample mean as MLE; this
led to the now famous \textit{Fisher--von
Mises--Langevin} (abbreviated FvML) distribution on spheres.

Of an entirely different nature is Campbell's \cite{C70} use of MLE
characterizations. In his
paper, Campbell establishes an equivalence between MLE
characterizations in the spirit of Gauss' principle and the minimum
discrimination information estimation of probabilities.
More recently, Puig \cite{P08} has
applied MLE characterizations in order to characterize the Harmonic
Law as the only statistical model on the positive real half-line that
satisfies a certain number of requirements. As a last example, we cite
the work of Ley and Paindaveine \cite{LP10a} who solved a long-standing
problem on skew-symmetric distributions through an argument resting on
Gauss' MLE characterization.

\subsection{Purpose of the paper} \label{secpurpose-paper}

As is perhaps intuitively clear, many of the results on MLE
characterizations stem from a common origin. As a matter of fact
many authors often follow the same ``smart path'' that can be
summarized in three steps: (a) choose the role of the parameter of
interest $\theta$ (location or scale); (b) choose a remarkable form
for the MLE for $\theta$ (e.g., sample mean, variance or median); (c)
use the freedom of choice in the samples as well as the sample size
(two samples of respective sizes 2 and 3, one sample of size 3, all
samples of all sizes$,\ldots$) to obtain the largest class of
distributions satisfying certain assumptions (continuous at a single
point, continuous, differentiable$,\ldots$) which share this specific
MLE. While similar, the arguments leading to the different results
nevertheless are largely ad hoc and rest upon crafty
manipulations of the explicit given form of the MLE. Moreover, step
(c) contains assumptions on the minimal sample size and on the
properties of the distributions being characterized, the necessity of
which is barely addressed.

The purpose of the present paper is to unify this important literature
by explaining the mechanism behind this ``smart path'' in the more
general context of one-parameter \textit{group families} (as defined in
Section \ref{secother-parameters}) although still with particular
focus on location and scale families. In doing so, we will introduce
the concept of {minimal covering sample size} ({MCSS}), a quantity
whose value depends on the structure of the support of the target
distribution and from which one deduces the a priori minimal
necessary sample size (MNSS) at which an MLE characterization holds
for a given family of distributions (see Duerinckx and Ley \cite{DL12}
where this notion
was introduced).
As we shall see, the MCSS and MNSS 
explain many of the differences in the
``always''-statements appearing throughout the literature on MLE
characterization. Moreover our unified perspective on MLE characterizations
not only permits to identify the minimal sufficient conditions under which
group families are characterized by their MLEs
(MLE-characterizable)
but also provide tools for (easily) constructing new MLE
characterizations of
many important families of distributions.

In a nutshell, our goal is to (i) propose a unified perspective on
MLE characterizations for one-parameter group families, (ii) answer
the question of which such families are MLE-characterizable and find
their MNSS, (iii) retrieve, improve on and better understand existing
results via our general analysis, and (iv) construct new MLE
characterization results.
Our contribution to the important literature on MLE characterizations
unifies the cornerstone references Teicher \cite{T61}, Ferguson \cite
{F62}, Kagan, Linnik and Rao \cite{KLR73}, Marshall and Olkin \cite
{MO93}, Azzalini and Genton \cite{AG07}
(to cite but these) and complements the understanding brought by the
seminal works of Poincar{{\'e}} \cite{P12}, Campbell \cite{C70} and
Bondesson \cite{B97}.

\subsection{Outline of the paper} \label{secoutline}

In Section \ref{notations} we describe the framework of our study, give
all necessary notations and introduce the so-called \textit{equivalence
classes}. In Section \ref{secminim-necess-sample} we establish and
interpret the above-mentioned notion of MCSS which will be central to
this paper. In Section \ref{secloc}, we derive the MLE characterization
for univariate location families, while in Section \ref{secscale} we
proceed in a similar way with univariate scale families. In Section
\ref{secother-parameters} we obtain MLE characterizations for general
one-parameter group families, allowing us to study other roles of the
parameter (e.g., skewness). In Section \ref{Appl} we apply our findings
to particular families of distributions. We conclude the paper with a
discussion, in Section \ref{secdiscu}, of the different possible
extensions that are yet to be explored.

\section{Notations and generalities on ML estimators, equivalence
classes}\label{notations}

Throughout we consider observations $\mathbf X^{(n)} =(X_1, X_2,\ldots,
X_n)$ that are sampled independently from a distribution $P_\theta^f$
(with density $f$) which we suppose to be entirely known up to a
parameter $\theta\in\Theta\subset\R$. The true parameter value
$\theta_0\in\Theta$ is estimated by ML estimation on basis of
$\mathbf{X}^{(n)}$. As explained in the \hyperref[intro]{Introduction},
our aim consists
in determining which classes of distributions are identifiable by means
of the MLE of the parameter of interest $\theta$, a parameter that can,
in principle, be of any nature (i.e., location, scale, etc.). On the
target family of distributions $\{P_\theta^f\dvtx  \theta\in\Theta\}$ we
make the following general assumptions:
\begin{enumerate}[-]
\item[-](A1) The parameter space $\Theta$ contains an open set
$\Theta_0$ of which the true parameter $\theta_0$ is an interior
point.
\item[-](A2) For each $\theta\in\Theta_0$, the distribution $P_\theta^f$
has support $S$ independent of $\theta$.
\item[-](A3) For all $1\le i \le n$ the random variable $X_i$ has a
density $f(x_i;\theta)$ with respect to the (dominating) Lebesgue
measure.
\item[-](A4) For $\theta\neq\theta' \in\Theta_0$, we have $P_\theta^f
\neq P_{\theta'}^f$.
\end{enumerate}
These assumptions are taken from Lehmann and Casella \cite{LC98}, page 444.

\begin{remark} \label{remark1}
Although typically we will be concerned
with either \textit{location families} with densities of the form
$f(x; \theta) = f(x-\theta)$ for $\theta\in\R$ the location
parameter, or \textit{scale families} with densities of the form $f(x;
\theta) = \theta f(x \theta)$ for $\theta\in\R^+_0$ the scale
parameter, other roles for $\theta$ (skewness, tail behavior\ldots)
can also be considered; see Section \ref{secother-parameters} where
we detail our approach for one-parameter group families.
\end{remark}
%
\begin{remark}
While throughout the paper we restrict our attention to the univariate setting,
it is easy to see that our arguments are in some cases transposable word-by-word
to the multivariate case. We discuss this matter briefly in
Section \ref{secdiscu}.
\end{remark}
%
\begin{remark} Assumption (A2) implies
that only densities with full support $\R$ may be considered for ML
estimation of a location parameter, while only densities with
support either $\R$, $\R^+$ or $\R^-$ may be considered for ML
estimation of a scale parameter. Despite the fact that these
restrictions are natural in the present context they can, if deemed
necessary, be lifted. We will briefly discuss this topic in
Section \ref{secdiscu}.
\end{remark}
We define, for a fixed sample size $n\ge1$, the MLE of the parameter
$\theta$ as (if it exists) the measurable function
\[
\hat\theta^{(n)}_f\dvtx  S^n:= S \times\cdots\times
S \to\Theta_0\dvtx  \mathbf x^{(n)}:= (x_1,\ldots,
x_n) \mapsto\hat\theta^{(n)}_f\bigl(\mathbf
x^{(n)}\bigr)
\]
for which
%
\begin{equation}
\label{eqmlenef} \prod_{i=1}^n f
\bigl(x_i; \hat\theta^{(n)}_f\bigl(\mathbf
x^{(n)}\bigr)\bigr) \ge\prod_{i=1}^n
f(x_i; \theta)
\end{equation}
for all $\theta\in\Theta_0$ and all samples $\mathbf x^{(n)} \in
S^n$ of size $n$.
It is not trivial to provide minimal conditions on $f$ under which
$\hat\theta^{(n)}_f(\mathbf x^{(n)})$ exists, is uniquely defined and
satisfies the necessary measurability conditions.
Consistency of the MLE is also a delicate matter and further regularity
conditions are required for the problem to make sense. As in
Cram{\'e}r \cite{C46,C46b} one may suppose that, for almost all $x$,
the density
$f(x; \theta)$ is differentiable with respect to $\theta$. This
allows to define the MLE as the solution of the local likelihood equation
%
\begin{equation}
\label{eqmleneflocal}\sum_{i=1}^n
\varphi_f(x_i; \theta) = 0,
\end{equation}
where
\[
\varphi_f(x; \theta):= \frac{\partial}{\partial\theta} \log f(x; \theta)
\]
is the \textit{score function} of the density $f$ associated with the
parameter $\theta$ (we set this function to 0 outside the support of
$f$). The solution to (\ref{eqmleneflocal}) has, at least
asymptotically, the required properties (see Lehmann and Casella \cite
{LC98}, Theorem
6.3.7). Furthermore, this way of proceeding allows for a simple
sufficient condition for uniqueness of the MLE: the mapping
$x\mapsto\varphi_f(x; \theta)$ has to be strictly monotone and to
cross the $x$-axis. Note
that this requirement coincides with \textit{strong unimodality} or
\textit{log-concavity} of the density $f$ when $\theta$ is a location
parameter (see Lehmann and Casella \cite{LC98}, Exercise 6.3.15).

\begin{remark}
Although in general there is no explicit expression for the MLE of a
given parametric family, there exist several important distributions
that not only satisfy all the above
requirements but also allow for MLEs which take on a remarkable
form. Taking $f = \phi$ the standard normal density,
the MLE for the location parameter is $\bar x:=
\frac1n\sum_{i=1}^nx_i$, the sample arithmetic mean, while that for
the scale parameter is $(\frac1n\sum_{i=1}^nx_i^2)^{1/2}$,
the square root of the sample arithmetic mean of squares. Taking $f$ the
exponential density, the MLE for the scale parameter becomes
$\bar{x}$.
\end{remark}

As outlined in Section \ref{secpurpose-paper}, our objective in this paper
is to determine minimal conditions under which a given form of MLE for
a given type of parameter identifies a specific probability
distribution. Upon making this statement there immediately arises a
trivial identification problem (due to our definition of maximum
likelihood estimators) which we first need to evacuate. Indeed if
$\hat{\theta}_f^{(n)}$ maximizes the $f$-likelihood function, then it
also maximizes the $g$-likelihood function for any function $g=cf^d$
with $d>0$ and $c$ a normalizing constant. In fact, from our
definitions (\ref{eqmlenef}) and (\ref{eqmleneflocal}), it is
immediate that any two parametric densities $f(x;\theta)$ and
$g(x;\theta)$ with same support $S$ and such that
%
\begin{equation}
\label{scoreeq} \varphi_{g}(x;\theta)=d \varphi_{f}(x;
\theta) \qquad\forall x\in S
\end{equation}
for some $d>0$ share the same MLE for $\theta$. This seemingly trivial
and innocuous observation leads us to the introduction of a concept
that happens to be fundamental for this paper and for MLE
characterizations in general: the concept of (parameter-dependent)
equivalence classes (e.c. hereafter), meaning that two densities
are equivalent if their score functions satisfy (\ref{scoreeq}). It is
obvious that e.c.'s constitute a partition of the space of
distributions. Stated simply, without further
conditions or specifications on the density associated with the target
distribution $P_\theta^f$ (the target density $f$) and on the
functions $g$, MLE characterization theorems
identify an e.c. of distributions rather than a single
well-specified distribution.

\begin{remark}
Considering, for example, Gauss' MLE characterization, both Teicher \cite{T61}
and Azzalini and Genton \cite{AG07}, to cite but these, identify the Gaussian
distribution with respect to its location parameter only up to an
unknown variance. On the contrary, when dealing with scale
characterizations, Teicher \cite{T61} imposes a tail-constraint in
order to
be able to single out the standard exponential and the standard
Gaussian distribution.
\end{remark}

The use of e.c.'s are far from new in the context of ML estimation.
Indeed Daniel Bernoulli, who is considered to be one of the first
authors to introduce the idea underlying the concept of maximum
likelihood, noted as early as in 1778 that the roots of his
``likelihood function'' would not change by squaring the semi-circular
density he used (see Stigler \cite{S99}, Chapter 16). Further ideas of
Bernoulli can be found in Kendall \cite{K61}, Section 19, and they clearly
demonstrate that the equivalence (\ref{scoreeq}) can be viewed as a
modern expression of Bernoulli's early thoughts.

In what follows, we shall state our results in the most general
possible way, without
(at least in the main results) considering additional identification
constraints. As will become clear from the subsequent sections, the
nature of the
parameter of interest $\theta$ heavily influences the partition of the
space of distributions, so
that, for each type of parameter, one needs to identify the
e.c.'s (see the beginnings of Sections \ref{secloc}
and \ref{secscale} for an illustration).

The framework we have developed so far allows us to reformulate the
question underpinning
the present article in a more transparent form, namely
``\textit{Do there exist two distinct e.c.'s $\mathcal{F}(\theta)$ and
$\mathcal{G}(\theta)$ such that the distributions $P_\theta^{f}$
and $P_\theta^{g}$ for $f\in\mathcal{F}(\theta)$ and
$g\in\mathcal{G}(\theta)$ share a given MLE of the parameter of
interest $\theta$}''?
As is perhaps already clear, the answer to this question is
negative~-- at least in the interesting cases. To see why this ought to
be the case
let $\mathcal{F}(\theta)$ and $\mathcal{G}(\theta)$
be two e.c.'s
which share an MLE for some sample size $n\ge1$, in
other words suppose that, for some $n\ge1$, some $f\in\mathcal
{F}(\theta)$ and
some $g\in\mathcal{G}(\theta)$, the estimator $\hat{\theta}^{(n)}_g$
coincides with $\hat{\theta}^{(n)}_f$, that is,
%
\begin{equation}
\label{eqmlenefg} \sum_{i=1}^n\log g
\bigl(x_i; \hat{\theta}^{(n)}_f\bigl(
\mathbf{x}^{(n)}\bigr)\bigr) \ge\sum_{i=1}^n
\log g(x_i; \theta)
\end{equation}
for all $\theta\in\Theta_0$ and all samples $\mathbf x^{(n)} \in
S^n$ of size $n$. Then clearly the only way for $f$ and $g$ to
satisfy (\ref{eqmlenefg}) for all samples of size $n$ is
that they be strongly related to one another; as we will
see, under reasonable conditions, one can go one step further and
deduce that if (\ref{eqmlenefg}) holds true for a ``sufficiently
large $n$'' then $f$ and $g$ must belong to the same
e.c. This intuition
is the heart of all the literature on MLE
characterizations.

\section{The minimal covering sample size}\label{secminim-necess-sample}

The first step towards establishing our general MLE characterizations
is to
gain a better understanding of the meaning of ``sufficiently large
$n$'' for (\ref{eqmlenefg}) to induce a
characterization theorem.
%
We use the notation $H_{n}$ to denote the hyperplane
\[
H_{n} = \Biggl\{\mathbf b^{(n)}=(b_1,\ldots,
b_{n}) \in\R^n \bigg| \sum_{j=1}^n
b_j = 0 \Biggr\}
\]
and associate, with any parametric distribution $P_\theta^f$
satisfying the
requirements of Section \ref{notations}, the collection(s) of sets
%
\begin{equation}
\label{eqsetsBn} B_n^{f, \theta}(A) = \bigl\{\mathbf
b^{(n)}\in H_n | \exists\mathbf x^{(n)} \in
A^n \mbox{ with } b_j = \varphi_f(x_j;
\theta) \mbox{ for all } 1\le j \le n \bigr\},
\end{equation}
where $A^n = A \times\cdots\times A$ is the $n$-fold cartesian
product of $A \subseteq S$, the support of the target $f$, and
$\theta\in\Theta_0$.
The interplay between the sets $B_n^{f, \theta}$ and the hyperplanes
$H_n$ determines the minimal sample size $n$ for a characterization
theorem to hold. The following lemma is crucial to our approach.



\begin{lemma} \label{lemmalem2} Suppose that for some
$\theta\in\Theta_0$ the mapping $x \mapsto\varphi_f(x; \theta)$ is
strictly monotone over some interval $\mathcal{X}\subset S$ and that
the (restricted) image $\mathcal{I}_{f, \theta}(\mathcal X):=
\{\varphi_f(x; \theta) | x \in\mathcal{X}\}$ is of the form
$(-P^{-}_{f, \theta},P^{+}_{f, \theta})$, for positive\vadjust{\goodbreak} constants
$P^{-}_{f, \theta}, P^{+}_{f, \theta}$ (possibly infinite). Then, for
all $n\ge1$,
$B_n^{f,\theta}(\mathcal{X})=H_{n}\cap
(\mathcal{I}_{f,\theta}(\mathcal X))^{n}$. Also, letting
%
\begin{equation}
\label{eqMNSS} N_{f,\theta}=\cases{ 2, &\quad if $P^{-}_{f, \theta}=
P^{+}_{f,\theta}$,
\vspace*{2pt}\cr
\displaystyle \biggl\lceil\frac{\max(P^+_{f, \theta},P^-_{f, \theta})}{\min(P^+_{f,
\theta},P^-_{f, \theta})}+1 \biggr
\rceil, &\quad if $P^{-}_{f, \theta}\neq P^{+}_{f, \theta}$
and $P^{-}_{f, \theta}, P^{+}_{f, \theta} < +
\infty$,
\vspace*{2pt}\cr
\infty, &\quad otherwise,}
\end{equation}
we get:
\begin{itemize}
\item for $n < N_{f,\theta}$, the orthogonal projections
$\Pi_{x_{j}}(B_n^{f,\theta}(\mathcal{X}))\subsetneq(-P^-_{f, \theta},
P^+_{f, \theta})$ for all $j=1,\ldots, n$;
\item for all $n \ge N_{f,\theta}$,
$\Pi_{x_j}(B_n^{f,\theta}(\mathcal{X})) =(-P^-_{f, \theta}, P^+_{f,
\theta})$
for all $j = 1,\ldots, n$.
\end{itemize}
The number $N_{f,\theta}$ is called the minimal covering sample size
(MCSS) associated to the interval~$\mathcal{X}$.
\end{lemma}

The MCSS, which is by definition greater than 2, is the
smallest possible value of the sample size $n$ that ensures that all
the projections of
$B_n^{f,\theta}(\mathcal{X})$ onto the distinct subspaces generated by
each observation $x_j, j=1,\ldots,n$, cover entirely
$\mathcal{I}_{f,\theta}(\mathcal X)$. For ease of reference we will
say, whenever this property holds true, that
$B_n^{f,\theta}(\mathcal{X})$ is \textit{projectable}.



\begin{example} \label{exfs-bn-loc}
Take $f(x;\theta) = \phi(x-\theta)$ the standard Gaussian density
with $\theta\in\R$ a location parameter. We have $S=\R$,
$\varphi_\phi(x;\theta) = x-\theta$ (and $\hat\theta_\phi^{(n)}=
\bar x$,
the sample mean). Then $\varphi_\phi(x;\theta)$ is invertible over
$\R$, $\mathcal{I}_{\phi,\theta}(\R)= \R$ (for all $\theta$) and easy
computations show that $B_n^{\phi,\theta}(\R)=H_n$ for all $n$. Note
that we always have
$\Pi_{x_{j}}(B_n^{\phi,\theta}(\R)) = \R$; in other words,
$B_n^{\phi,\theta}(\R)$ is projectable for all $n\ge2$ (and all
$\theta$) and hence $\mbox{MCSS}=2$.
\end{example}

\begin{example} \label{exfs-bn-toy}
Take $\varphi_{f}(\cdot;\theta)$ to be monotone on $S=\R$ with
symmetric image $(-1,1)$, say, independently of $\theta$. Then clearly
$B_2^{f, \theta}(\R)$ is the intersection between the line $H_2\equiv
x+y=0$ (see Figure \ref{Fig1}) and the square $(-1,
1)^2$, while $B_3^f(\R)$
is the intersection between the plane
$H_3\equiv x+y+z=0$ and the cube $(-1, 1)^3$.
%
\begin{figure}

\includegraphics{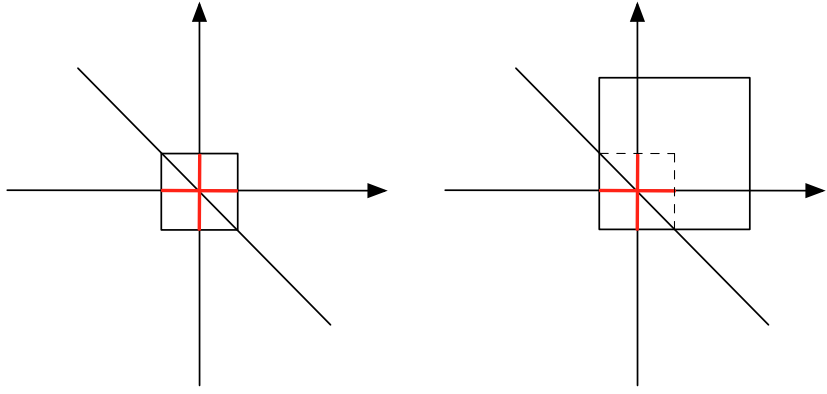}

\caption{The sets $B_{2}^{f, \theta}(\R)$ when $\mathcal{I}_{f,
\theta} =[-1, 1]$
(left plot) and $\mathcal{I}_{f, \theta}=[-1, 3]$ (right plot). The
red lines
indicate the marginal projections.}
\label{Fig1}
\end{figure}
Consequently, the coordinates of points on $B_2^{f, \theta}(\R)$ and
$B_3^{f, \theta}(\R)$ cover the full image $(-1,1)$. Hence, $\mbox{MCSS}=2$.
\end{example}

\begin{example}
Take $\varphi_f(\cdot;\theta)$ to be monotone on $S =
\R$ with skewed image
$(-1,3)$, say, independently of $\theta$. Then, while
$B_2^{f, \theta}(\R)$ and $B_3^{f, \theta}(\R)$ remain defined as in
Example \ref{exfs-bn-toy} (with $(-1,3)$
replacing $(-1,1)$), coordinates of points in these domains do not
cover the full image. In fact, in $B_2^{f, \theta}(\R)$ these
coordinates only
cover the interval $(-1,1)$ (see Figure~\ref{Fig1}), in $B_3^{f,
\theta}(\R)$
these coordinates only cover the interval $(-1, 2)$ and it is only
from $n\ge4$ onwards that the coordinates of points in $B_n^{f, \theta
}(\R)$
cover the full image. Hence, $\mbox{MCSS}=4$.
\end{example}

\begin{example} \label{exfs-bn-scale}
Take $f(x;\theta) = \theta\phi(\theta x)$ the Gaussian density with
$\theta\in\R^+_0$ a scale parameter. We have $S=\R$,
$\varphi_\phi(x;\theta) = \frac1\theta(1-\theta^2x^2)$. Then
$\varphi_\phi(\cdot;\theta)$ is invertible over $\R_0^+$ and $\R_0^-$,
separately, and $\mathcal I_{f,\theta}(\R_0^{\pm})=(-\infty,1/\theta)$.
Note that we have $\Pi_{x_{j}}(B_n^{f,\theta}(\R_0^{\pm})) =
(-(n-1)/\theta, 1/\theta)$ for all $j$, all $\theta$ and all $n\ge2$.
In other words, $B_n^{f,\theta}(\R_{0}^{\pm})$ is only asymptotically
projectable. Hence, $\mbox{MCSS} = + \infty$.
\end{example}

\begin{pf*}{Proof of Lemma \ref{lemmalem2}}
The assumptions on the image of $\varphi_f(\cdot;\theta)$ guarantee
the existence of a point in $\mathcal{X}$ where $\varphi_f$ crosses
the $x$-axis so that $B_1^{f,\theta}(\mathcal{X}) =
\{0\}$. Also, for all $n\ge1$, we have
$B_n^{f,\theta}(\mathcal{X})=H_{n}\cap
(\mathcal{I}_{f,\theta}(\mathcal X))^{n}$; this follows by
definition. Regarding the MCSS, first take $P^{-}_{f, \theta
}=P^{+}_{f,\theta} = P$ (possibly
infinite). Then, for all $n\ge2$, $H_n\cap(-P, P)^n$ contains for
each of the $n$ coordinates the full interval $(-P, P)$ (see Figure
\ref{Fig1}); hence $\mbox{MCSS} = 2$.
Next suppose that $P^{-}_{f, \theta} < P^{+}_{f,\theta} <\infty$ and
consider a
``worst-case scenario'' by taking a point at the extreme of
$H_n\cap(\mathcal{I}_{f,\theta}(\mathcal X))^{n}$, with one
coordinate set to $b_1 = P^+_{f, \theta}-\varepsilon$ for some $\varepsilon
>0$. Then,
in order to construct a sample satisfying $\sum_{i=1}^n b_i = 0$, it is
necessary to choose the remaining $n-1$-tuple $(b_2,\ldots,b_n)$ so
as to satisfy $\sum_{i=2}^nb_i=-b_1$. Since the best choice in this
respect consists in setting all $b_i$ near the other extremum
$-P^-_{f, \theta}+\varepsilon'$ for $\varepsilon'>0$, we see that,
depending on
the magnitude of the ratio $P^+_{f, \theta}/P^-_{f, \theta}$, a given
sample size $n$ may
not\vspace*{1pt} be large enough for the equality to hold. In order to
palliate this it suffices to take $N_{f,\theta}$ to be the smallest natural
number such that $P^+_{f, \theta}-(N_{f,\theta}-1)P^-_{f, \theta
}\leq0$, that is,
$N_{f,\theta}= \lceil{P^+_{f, \theta}}/{P^-_{f,
\theta}}+1 \rceil$. The case $P^+_{f, \theta}< P^-_{f, \theta
}<\infty$
follows along the same lines, and hence
\[
\mbox{MCSS} = \biggl\lceil\frac{\max(P^+_{f, \theta},P^-_{f, \theta
})}{\min(P^+_{f,
\theta},P^-_{f, \theta})}+1 \biggr\rceil.
\]
The same argumentation applies in the case where either one of
$P^{+}_{f, \theta}$ or $P^{-}_{f, \theta}$ is infinite, this time with
$\mbox{MCSS} = +\infty$. This concludes the proof.
\end{pf*}

Now look at the connection between the MCSS and MLE
characterizations. Let $f$ and $g$ be two representatives of distinct
e.c.'s with $f$ the target density. Under the assumption of
$\theta$-differentiability of $g$, the defining equation (\ref
{eqmlenefg}) can
be re-expressed as
\[
\sum_{i=1}^n\varphi_g
\bigl(x_i;\hat{\theta}^{(n)}_f\bigl(\mathbf
x^{(n)}\bigr)\bigr)=0\qquad\mbox{for all } \mathbf x^{(n)}\in
S^n,
\]
%
which in turn can be rewritten as
%
\begin{equation}
\label{eqs1} \sum_{i=1}^n
\varphi_g(x_i;\theta)=0 \qquad\mbox{for all } \theta\mbox{
and all } \mathbf{x}^{(n)}\mbox{ such that } \sum
_{i=1}^n \varphi_f(x_i;
\theta)= 0
\end{equation}
(here $\theta$ and $\mathbf{x}^{(n)}$ are interdependent) or, equivalently,
%
\begin{equation}
\label{MLEeq} \sum_{i=1}^n
h(y_{i};\theta)=0 \qquad\mbox{for all } \theta\mbox{ and all }
\mathbf{y}^{(n)} \in B_n^{f,\theta}(\mathcal{X})
\end{equation}
for each interval $\mathcal{X}$ on which $y\mapsto\varphi_f(y;\theta)$
is invertible for all $\theta$, with $h(y;\theta) =
\varphi_g(\varphi_f^{-1}(y;\allowbreak \theta); \theta)$. 
Equation (\ref{MLEeq}) completely identifies the function $h$ (and
hence the interconnection between $g$ and $f$), at least when
$B_n^{f,\theta}(\mathcal{X})$ is sufficiently
rich. This richness depends strongly on the MCSS introduced in Lemma
\ref{lemmalem2}. Indeed supposing that (\ref{MLEeq}) is only valid
for a sample size smaller than the MCSS implies that portions of the
images $\mathcal{I}_{f,\theta}(\mathcal{X})$ cannot be
reached, so that $h$ cannot be identified over its entire
support. This necessarily implies that $\mbox{MNSS} \ge
\mbox{MCSS}$. It is, however, pointless to try to solve
(\ref{MLEeq}) in all generality and it is now necessary to
specify the role of $\theta$ in order to pursue.
%
%
%
We will do so in detail in the next sections, first in the case of
location parameters; our arguments will afterwards adapt directly to
other parameter choices.

\section{MLE characterization for location parameter
families}\label{secloc}

We start by identifying the e.c.'s for $\theta$ a location
parameter. In such a case, $\Theta_0= S=\R$, and the location score
functions are of the form
$\varphi_f(x;\theta)=\varphi_f(x-\theta)=-f'(x-\theta)/f(x-\theta)$
over $\R$, so that equation (\ref{scoreeq}) turns into a simple
first-order differential equation whose solution yields $g(x)=c
(f(x))^d$ for some $d>0$ and $c$ the normalizing constant. Thus, all
densities which are linked one to another via that relationship belong
to a same e.c. We here attract the reader's attention to the fact
that, for $f=\phi$ the standard Gaussian density, such transformations
reduce to a non-specification of the variance, which is clearly in
line with Gauss' MLE characterization as stated by Teicher \cite{T61}
or Azzalini and Genton \cite{AG07}.

Our first main theorem is, in essence, a generalization of Gauss'
MLE characterization from the Gaussian distribution to the entire class
of log-concave
distributions with continuous score function.

\begin{theorem}\label{Theoloc}
Let $\mathcal{F}(\mathrm{loc})$ and $\mathcal{G}(\mathrm{loc})$ be two
distinct location-based e.c.'s and let their respective representatives
$f$ and $g$ be two continuously differentiable densities with full
support $\R$. Let $x\mapsto\varphi_f(x) = -f'(x)/f(x)$ be the location
score function of $f$. If $\varphi_f$ is invertible over $\R$ and
crosses the $x$-axis then there exists $N \in\N$ such that, for any
$n\ge N$, we have $\hat{\theta}^{(n)}_f=\hat{\theta}^{(n)}_g$ for all
samples of size $n$ if and only if there exist constants $c,d\in\R_0^+$
such that $g(x)=c(f(x))^d$ for all $x\in\R$, that is, if and only if
$\mathcal{F}(\mathrm{loc})=\mathcal{G}(\mathrm{loc})$.
The smallest integer for which this holds (the
minimal necessary sample size) is
$\mbox{MNSS}= \max\{N_f, 3\}$,
with $N_f$ the MCSS as defined in Lemma~\ref{lemmalem2}.
\end{theorem}

\begin{pf}
The sufficient condition is trivial. To prove necessity first note
how our assumptions on $f$ ensure that the score function $\varphi_f$
is strictly increasing on the whole real line $\R$ and has a unique
root. This allows us
to write the image $\operatorname{Im}(\varphi_f)$ as $(-P^-_f,P^+_f)$ with
$0<P^-_f,P^+_f\leq\infty$.
The differentiability of $g$ and the nature of the parameter $\theta$
permit us to rewrite, for any admissible $\theta$, (\ref{eqs1}) as
%
\begin{equation}
\label{MLEeqloc} \sum_{i=1}^n
\varphi_g(x_i - \theta)=0 \qquad\mbox{for all } \mathbf
{x}^{(n)}\in\R^n \mbox{ such that } \sum
_{i=1}^n \varphi_f(x_i-
\theta)= 0,
\end{equation}
where $\varphi_g(x)=-g'(x)/g(x)$.
Using the strict monotonicity of $\varphi_f$ one then concludes
that (\ref{MLEeqloc}) is equivalent to requiring that $g$ satisfies
%
\begin{equation}
\label{Starteq} \sum_{i=1}^nh(b_i)=0
\qquad\mbox{for all } (b_1,\ldots,b_n)\in B_n^{f, \theta}(
\R),
\end{equation}
where $h=\varphi_g\circ\varphi_f^{-1}$ and $b_i=\varphi
_f(x_i-\theta)$, $i=1,\ldots,n$, as
in (\ref{eqsetsBn}). In what follows, we shall use our liberty of
choice among all $n$-tuples $\mathbf{b}^{(n)}\in B_n^{f, \theta}(\R
)$ in order to gain sufficient information on $h$ to conclude.

First,\vspace*{1pt} suppose that $P^-_{f, \theta}=P^+_{f, \theta}=P$, hence that
the image of $\varphi_f$ is symmetric. We know from
Lemma \ref{lemmalem2} that the corresponding MCSS equals 2, hence
that two observations suffice to make $B_n^{f, \theta}(\R)$
projectable. Therefore, for any $n\ge2$, we can always
build an $n$-tuple $b_1,\ldots,b_n$ such that $b_2=-b_1$ for all $b_1
\in(-P, P)$ and $b_i=0$ for $i=3,\ldots,n$. From (\ref{Starteq})
we then deduce that $h$ satisfies the equality $h(-a)=-h(a)$ for all
$a\in(-P,P)$, hence that $h$ is odd on $(-P,
P)$. Evidently this leaves $h$ undetermined, hence the MNSS must at
least equal 3. For
$n\geq3$, choose an $n$-tuple such that
$b_3=-b_1-b_2$ and $b_i=0$ for $i=4,\ldots,n$, for
$b_1,b_2\in(-P,P)$ such that $b_1+b_2\in(-P,P)$. Using
(\ref{Starteq}) combined
with the antisymmetry of $h$ we deduce that this function must
satisfy
%
\begin{equation}
\label{Cauchy} h(b)+h(c)=h(b+c)
\end{equation}
for all $b,c \in(-P,P)$ such that $b+c\in(-P,P)$. One recognizes in
(\ref{Cauchy}) a (restricted) form of the celebrated \textit{Cauchy
functional equation}. Assume that $P<\infty$; then
$h(P/2)$, say, is finite and standard
arguments (see, e.g.,
Acz{\'e}l and Dhombres \cite{AD89}), imply that our solution $h$
satisfies $h(uP/2)=uh(P/2)$
for all $u\in(-2,2)$ and we conclude that $h(x)= d x$ for
all $x\in(-P,P)$, with $d = h(P/2)/(P/2)\in\R$. Considering
$x=\varphi_f(y)$ for $y\in\varphi_f^{-1}(-P,P)=\R$, we obtain that
$\varphi_g(y)=d \varphi_f(y)$. Solving
this first-order differential equation gives $g(y)=c(f(y))^d$ for
all $y\in\R$, with $c$ a constant. In order for the function $g$
to be integrable over $\R$, the constant $d$ must be strictly
positive; in order for $g$ to be positive and integrate to 1, the
constant $c$ must be a
normalizing constant. Thus, for
$P<\infty$, the problem is solved. For $P=\infty$, the situation
becomes even simpler as (\ref{Cauchy}) is then \textit{precisely} the Cauchy
functional equation, and one may immediately draw the same
conclusion as for finite $P$.

Let us\vspace*{1pt} now consider the case where $\varphi_f$ has a skewed
image and set $P=\min(P^-_f,P^+_f)$ (note that $P$ is necessarily
finite as otherwise
$\operatorname{Im}(\varphi_f)$ would be symmetric). First restricting our
attention to $(-P,P)$, we can repeat the above
arguments to deduce that $g(y)=c(f(y))^d$ for
all $y\in\varphi_f^{-1}(-P,P)\subsetneq\R$. We thus further need to
investigate the behavior of $h$ on the remaining part of $\operatorname{Im}(\varphi_f)$ which, for the sake of simplicity, we denote as
$\operatorname{Out}(P)$ (it is either $(-P^-_f,-P)$ or $(P,P^+_f)$). To this end,
we precisely need to know the MCSS and hence call upon
Lemma \ref{lemmalem2}. Fixing $n\geq N_f$
and taking a sample
$(b_1,\ldots,b_n)$ such that $\sum_{i=1}^nb_i=0$ with
$b_1\in\operatorname{Out}(P)$ and $(b_2,\ldots,b_n)\in(-P,P)^{n-1}$, we can
apply (\ref{Starteq}) to get $h(b_1)+\sum_{i=2}^nh(b_i)=0$ and
hence, from our knowledge about the behavior of $h$ on $(-P, P)$,
we deduce that
\[
h(b_1)=-\sum_{i=2}^nb_ih(P)/P=b_1h(P)/P,
\]
since $h(P)$ is necessarily finite. Consequently we get $h(y)=yh(P)/P$
for all $y\in(-P,P)\cup\operatorname{Out}(P)=\operatorname{Im}(\varphi_f)$ and $g(y)=c(f(y))^d$ for
all $y\in\R$, and the conclusion follows.

The proof of the theorem is nearly complete: all that remains is to
show that the $\mathrm{MNSS}=\max\{3,N_f\}$ is minimal and sufficient.
The latter is immediate since if the result holds true for any sample
of size $N = \mathrm{MNSS}$ then, for any larger sample size $M >
\mathrm{MNSS}$, one can always consider $ \mathbf{x}^{(M)}$ such that
$\varphi_f(\mathbf{x}^{(M)}) = \mathbf{b}^{(M)} \in B_M^{f,
\theta}(\R)$ is of the form $(b_1,\ldots, b_N, 0,\ldots, 0)$ and
$(b_1,\ldots, b_N) \in B_N^{f, \theta}(\R)$, and work as above to
characterize the density. To prove the minimality of the
$\mathrm{MNSS}$ it suffices to exhibit specific counter-examples. This
is done in Examples \ref{ex1} and~\ref{ex2} below.
\end{pf}

%

%
\begin{example}\label{ex1} To see that $N=3$ is minimal when $\operatorname{Im}(\varphi_f)$ is symmetric, we need to construct two
distributions $g_{1}$, $g_2$ which share $f$'s MLE for all samples of size
2. Construct $g_1$ as in the proof of Theorem \ref{Theoloc}. To
construct $g_2$, it suffices to replace the function $h$ from
(\ref{Starteq}) with any odd function and to solve the resulting
equation in $g$ (while ensuring integrability of $g$).\vspace*{1pt} If, for
example, we choose $h(x)=d x^3$, then we readily obtain
$g(y)=c\exp(-d\int_{-\infty}^x(\varphi_f(y))^3\mrmd y)$; this is however
not a density for all $f$, though a good choice for $f=\phi$ the
Gaussian (for which $\varphi_\phi(x)=x$). Another way of
proceeding is to work as in Azzalini and Genton \cite{AG07} and to choose
$h(y)=y+w'(y)$ for some
differentiable even function $w$.
\end{example}

\begin{example}\label{ex2} Suppose that $\hat{\theta}^{(n)}_f=\hat
{\theta}^{(n)}_g$
for all samples of size $n$ for some $n<N_f$ when $N_f>3$. Then,
as is clear from Lemma \ref{lemmalem2}, the whole domain (i.e., $\R
$) of $f$ is not
identified by our technique and it suffices to choose any density which
is equal to $f$ on the maximal identifiable subdomain but differs
elsewhere. Expressed in terms of $h$ for the case $P^{-}_f<P^{+}_f$, we
can only identify $h$ on
some interval $(-P^-_f,P^-_f+a(n))$, say, with
$0<a(n)<P^+_f-P^-_f$. On the remaining part $(P^-_f+a(n),P^+_f)$, $h$ is
undetermined and hence can take any possible form, implying that the
relationship between $g$ and $f$ can only be established on the part
$\varphi_f^{-1}(-P^-_f,P^-_f+a(n))\subsetneq\R$.
\end{example}

As in Azzalini and Genton \cite{AG07}, it is sufficient to require in Theorem
\ref{Theoloc} that $g$ be continuously differentiable at a single
point for everything to run smoothly. Pursuing in this vein, it is
of course natural to enquire whether the result still holds if
no such regularity assumption is imposed on~$g$, that is, if we only
suppose that
the target density $f$ is differentiable but $g$ is a priori
not. Put simply the
question becomes that of enquiring whether the condition
%
\begin{equation}
\label{eqmlenefgloc} \sum_{i=1}^n\log g
\bigl(x_i - \hat\theta^{(n)}_f\bigl(
\mathbf{x}^{(n)}\bigr)\bigr) \ge\sum_{i=1}^n
\log g(x_i- \theta)
\end{equation}
for all $\mathbf{x}^{(n)}\in\R^n$ and all $\theta\in\R$ suffices
to determine
$g$. This is the approach
adopted, for example, in Teicher \cite{T61}, Kagan, Linnik and Rao \cite
{KLR73} or Marshall and Olkin \cite{MO93},
where it is shown
that having the likelihood condition (\ref{eqmlenefgloc}) with
$\hat\theta^{(n)}_f$ the sample mean
implies $g$ is the Gaussian as soon as the result holds for all
samples of sizes 2 and 3 simultaneously. Interestingly, in our
framework, this arguably more general
assumption on $g$ comes with a cost: our method of
proof then necessitates imposing
more restrictive assumptions on $f$ and requiring the likelihood
equations to hold for two sample sizes
\textit{simultaneously}.\looseness=1

\begin{theorem}\label{Theoloc2}
Let $\mathcal{F}(\mathrm{loc})$ and $\mathcal{G}(\mathrm{loc})$ be two
distinct location-based e.c.'s and let their respective
representatives $f$ and $g$ be two continuous densities with
full support $\R$. Suppose that $f$ is symmetric and continuously
differentiable, and assume that its location score function
$\varphi_f(x)$ is invertible over $\R$ and crosses the
$x$-axis. Then we have
$\hat{\theta}^{(n)}_f=\hat{\theta}^{(n)}_g$ for all samples of
sizes 2 and
$n\geq3$ simultaneously if and only if
there exist constants $c,d\in\R_0^+$ such that $g(x)=c(f(x))^d$ for
all $x\in\R$, that is, if and only if
$\mathcal{F}(\mathrm{loc})=\mathcal{G}(\mathrm{loc})$.\looseness=1
\end{theorem}

\begin{pf}
Our proof, which extends that of Teicher \cite{T61} from the Gaussian
case to
the entire class of symmetric log-concave densities $f$, proceeds in
two main steps: we first show that our assumptions on $g$ in fact
entail that $g$ is continuously differentiable, and then conclude by
applying Theorem \ref{Theoloc}. The additional sample size $n=2$
needed here stems from the first step.

Condition (\ref{eqmlenefgloc}) can be rewritten as
\[
\sum_{i=1}^n\log g(y_i) \ge
\sum_{i=1}^n\log g(y_i-
\theta)
\]
for all $\theta\in\R$ and $y_1,\ldots,y_n$ satisfying $\sum
_{i=1}^n\varphi_f(y_i)=0$. The latter expression in turn is equivalent to
%
\begin{eqnarray}
\label{Teich1} &&\sum_{i=1}^{n-1}\log
g(y_i)+\log g \Biggl(\varphi_f^{-1} \Biggl(-
\sum_{j=1}^{n-1}\varphi_f(y_j)
\Biggr) \Biggr)
\nonumber\\[-8pt]\\[-8pt]
&&\quad \geq\sum_{i=1}^{n-1}\log
g(y_i-\theta)+\log g \Biggl(\varphi_f^{-1}
\Biggl(-\sum_{j=1}^{n-1}\varphi_f(y_j)
\Biggr)-\theta\Biggr).\nonumber
\end{eqnarray}
Arguing as in Teicher \cite{T61}, it is sensible to confine our
attention at
first to symmetric densities~$g$. Using the assumed symmetric nature of
$f$ (and hence the oddness of $\varphi_f^{-1}$), considering the
sample size $n=2$ and setting observation $y_1$ equal to some $y\in\R
$, (\ref{Teich1}) simplifies into
%
\begin{equation}
\label{Teich2} 2\log g(y)\geq\log g(y-\theta)+\log g(y+\theta)
\end{equation}
for all $y,\theta\in\R$. Since $\log g$ is everywhere finite,
concave according to (\ref{Teich2}), and inherits measurability from
$g$, it is an \textit{a.e.}-continuously differentiable function.
Arrived at this point, we may apply Theorem \ref{Theoloc} to conclude
(note that the oddness of $\varphi_f$ makes Theorem \ref{Theoloc}
hold with MNSS equal to $3$).

Finally, for non-necessarily symmetric densities $g$, we can follow
exactly the argumentation from Teicher \cite{T61} and derive that the
previously obtained solution is the only one, hence the claim holds.
\end{pf}

We stress the fact that, as in Teicher \cite{T61}, we may further
weaken our
assumptions on $g$ by only requiring that it is lower semi-continuous
at the origin and need not have full support $\R$. Indeed, as shown in
Teicher's proof, continuity and a.e.-positivity ensue from the
above arguments.

One may wonder whether the symmetry assumption on the target density
$f$ is necessary or whether this second general location MLE
characterization theorem may in fact hold for the entire class of
log-concave densities as well. Our method of proof indeed
requires this assumption so as to enable us to deal with such
quantities as
$\varphi_f^{-1}(-\sum_{j=1}^{n-1}\varphi_f(y_j))$ in (\ref{Teich1});
without\vspace*{1pt} any assumption on $f$, for $n=2$, this expression does not
simplify into the agreeable form $-y_1$. Similarly, one may wonder
whether it is
necessary to suppose the result to hold for two sample sizes
simultaneously or whether one single sample size $N \ge3$ might not
suffice. We leave as open problems the question whether these
assumptions are necessary or simply
sufficient.

Finally, our Theorems \ref{Theoloc} and \ref{Theoloc2} do not cover
target densities whose location score function is monotone but not
invertible over the entire real line, that is, piecewise
constant. We do not consider explicitly such setups here since they do
not assure that the MLE is defined in a unique way. The strategy we
however suggest consists in applying
our results on the monotonicity intervals, to draw the necessary
conclusions and express $g$ in terms of $f$ on those intervals. If we
add the condition of monotonicity of $\varphi_g$, the equality
$\varphi_g(x)=d \varphi_f(x)$ has to hold over the entire support $\R$
as monotonicity imposes $\varphi_g$ to be constant outside the\vadjust{\goodbreak}
above-mentioned intervals. Since we here do not implicitly use the
intervals where $\varphi_f$ is constant, there might exist better
strategies, and consequently the smallest possible sample size we
obtain by following this scheme is an upper bound for the true
MNSS. The most extreme situation takes place when the target is a
Laplace distribution, in which case $\varphi_f(x)=\operatorname{sign}(x)$; we
refer to Kagan, Linnik and Rao \cite{KLR73} for a treatment of this particular
distribution. We will return to these matters briefly in Section \ref
{secdiscu}.

\section{MLE characterization for scale parameter families}\label{secscale}

As for location parameter families, we start by identifying the e.c.'s
when $\theta$ plays the role of a scale parameter. In such a
setup, $\Theta_0=\R_0^+$, $S=\R$, $\R^+_0$ or $\R^-_0$ in view of
assumption (A2), and the scale score functions are of the form
$\varphi_f(x;\theta)=\frac{1}{\theta}\psi_f(\theta x):=\frac
{1}{\theta}(1+\theta x
f'(\theta x)/f(\theta x))$ over $S$, so that equation (\ref{scoreeq})
turns into another quite simple first-order differential equation
whose solution leads to $g(x)=c |x|^{d-1} (f(x))^d$ for some $d>0$
(such that $g$ is integrable) and $c$ the normalizing constant. This
relationship defines the scale-based e.c.'s. It is to be noted that
$c=d=1$ when the origin belongs to the support, that is, when $S=\R$,
in which case the e.c.'s reduce to singletons~$\{f\}$.

Our main scale MLE characterization theorem is the exact equivalent of
Theorem \ref{Theoloc} with $\psi_f$ replacing $\varphi_f$, hence its
proof is omitted.

\begin{theorem}\label{Theosca}
Let $\mathcal{F}(\mathrm{sca})$ and $\mathcal{G}(\mathrm{sca})$ be two
distinct scale-based e.c.'s and let their respective representatives
$f$ and $g$ be two continuously differentiable densities with common
support $S$ (either $\R,\R_0^+$ or $\R_0^-$). Let $\psi_f(x) =
1+xf'(x)/f(x)$ be the scale score
function of $f$. If $\psi_f$ is invertible over
$S$ and crosses the $x$-axis then there exists
$N \in\N$ such that, for any $n\ge N$, we have
$\hat{\theta}^{(n)}_f=\hat{\theta}^{(n)}_g$ for all samples of size
$n$ if and only if
there exist constants $c,d\in\R_0^+$ such that
$g(x)=c|x|^{d-1}(f(x))^d$ for
all $x\in S$ (with $c=d=1$ for $S=\R$), that is, if and only if
$\mathcal{F}(\mathrm{sca})=\mathcal{G}(\mathrm{sca})$.
The smallest integer for which this holds is
$\mbox{MNSS}= \max\{N_f, 3\}$,
with $N_f$ the MCSS as defined in Lemma \ref{lemmalem2}.
\end{theorem}

As in the case of a location parameter,
requiring differentiability of the $g$'s is not
indispensable. One could indeed restrict the class of target
distributions under consideration, as in Teicher \cite{T61}. We leave
this as an easy exercise.

When dealing with scale families it is natural to work as in
Teicher \cite{T61} and add a scale-identification condition of the form
%
\begin{equation}
\label{eqsc-cond} \lim_{x\rightarrow0} {g(\lambda x)}/{g(x)}=\lim
_{x\rightarrow0} {f(\lambda x)}/{f(x)} \qquad\forall\lambda>0.
\end{equation}
Imposing this condition in Theorem
\ref{Theosca} allows, at least when the limit is finite, positive and
does not equal $1/\lambda$ (which occurs only for pathological cases;
we leave as an exercise to the reader to see why this type of limiting
behavior precludes identifications), to
deduce that $c=d=1$ for $S=\R_0^+$ and $\R_0^-$, and
hence $g=f$ in all cases. Interestingly Teicher already remarks that this
``\textit{seemingly ad hoc condition appears to be crucial}''; this is
clearly the case for a complete identification of the family of
densities which share a scale MLE.\vadjust{\goodbreak}

The invertibility condition imposed on $\psi_f$ is as natural in a
scale family context as the invertibility condition on $\varphi_f$ in
a location family setup (see Lehmann and Casella \cite{LC98}, page
502). Unfortunately it
suffers from one major drawback for $S=\R$: requiring invertibility of
$\psi_f$ over the whole real line forces us to discard several
interesting cases such as, for example, the standard normal density $\phi$,
for which
$\psi_\phi(x)=1-x^2$ is only invertible over the positive and
negative real half-lines, respectively. More generally any symmetric
density $f$ for
which $\varphi_f$ is invertible over $\R$ will suffer from that same
problem and hence will not be characterizable by means of
Theorem \ref{Theosca}. This flaw is nevertheless easily fixed,
since Lemma \ref{lemmalem2} is applicable even if $\psi_f$ is only
invertible over portions of its support. This leads to our next
general result (whose proof is omitted).

\begin{theorem}\label{Theosca2}
Let $\mathcal{F}(\mathrm{sca})$ and $\mathcal{G}(\mathrm{sca})$ be two
distinct scale-based e.c.'s and let their respective representatives
$f$ and $g$ be two continuously differentiable densities with full
support $\R$. Let the scale score function $\psi_f(x) = 1+xf'(x)/f(x)$
be invertible and cross the $x$-axis over $\R_0^+$ and~$\R_0^-$,
respectively. Then there\vspace*{1pt} exists $N \in\N$ such that, for any $n\ge N$,
we have $\hat{\theta}^{(n)}_f=\hat{\theta}^{(n)}_g$ for all samples of
size $n$ if and only if $g(x)=f(x)$ for all $x\in\R$.
Moreover the MNSS is given by $\max\{\mbox{MNSS}_-,\mbox{MNSS}_+\}$,
where $\mbox{MNSS}_-$ and $\mbox{MNSS}_+$, respectively, stand for the
MNSS required on each half-line.
\end{theorem}

It should be noted that the scale
condition (\ref{eqsc-cond}) is not necessary here since we are working
on the entire support $S=\R$ which imposes that $d=1$ as otherwise the
non-vanishing density $g$ would vanish at 0.

Finally note that the separation of the two real half-lines is tailored for
scale families because both $\R_0^+$ and $\R_0^-$ are invariant under the
action of the scale parameter, which permits us to work on each
half-line separately and put the ends together by continuity. The same
would not hold true for location families due to a lack of invariance,
that is, we could not ``glue
together'' location characterizations valid on complementary subsets
of the support.

\section{MLE characterization for one-parameter group families} \label
{secother-parameters}

The relevance of our approach is not confined to location and scale
families, but can be used for other $\theta$-parameter families with
$\theta$ neither a location nor a scale parameter. In this section, we
shall consider general one-parameter group families and provide them
with MLE characterization results. Group families play a central role
in statistics as they contain several well-known parametric families
(location, scale, several types of skew distributions as shown in Ley
and Paindaveine \cite{LP10b}$,\ldots$) and allow for significant
simplifications of the data under investigation (see Lehmann and
Casella \cite{LC98}, Section 1.4, for more details). To the best of the
authors' knowledge, there exist no MLE characterizations for group
families other than the location and scale families.

A univariate group family of distributions is obtained by subjecting a
scalar random variable with a fixed distribution to a suitable family
of transformations. More prosaically, let $X$ be a random variable
with density $f$ defined on its support $S$ and consider a
\textit{transformation group}~$\mathcal{H}$ (meaning that it is closed
under both composition and inversion) of monotone increasing functions
$H_\theta\dvtx D\subseteq\R\rightarrow S$ depending on a single real
parameter $\theta\in\Theta_0$. The family of random variables
$ \{ H_\theta^{-1}(X), H_{\theta} \in\mathcal{H} \}$ is
called a \textit{group family}. These variables possess
densities of the form
%
\begin{equation}
\label{group} f_\mathcal{H}(x;\theta):=H_\theta'(x)f
\bigl(H_\theta(x)\bigr),
\end{equation}
where $H_\theta'$ stands for the derivative of the mapping $x\mapsto
H_\theta(x)$ (which we therefore also suppose everywhere
differentiable); their support $D$ does not depend on $\theta$. We
call $\theta$ a $\mathcal{H}$-parameter for $f(x;\theta)$. The most
prominent examples are of course $\mathcal{H}_{\mathrm{loc}}:=\{H_\theta(x)=x-\theta, x,\theta\in\R\}$, leading to location
families, and $\mathcal{H}_{\mathrm{sca}}:=\{H_\theta(x)=\theta x,x\in
S,\theta\in\R_0^+\}$ for $S=\R,\R_0^+$ and $\R_0^-$, yielding scale
families. For further examples, we refer to Lehmann and Casella \cite
{LC98}, Section 1.4,
and the references therein.

Let us now determine the e.c.'s for $\mathcal{H}$-parameter
families. Assuming that the mappings
$\theta\mapsto H_\theta(x)$ and $\theta\mapsto H_\theta'(x)$ are
differentiable, the $\mathcal{H}$-score function associated with
densities of the form (\ref{group}) corresponds to
%
\begin{equation}
\label{H-score} \varphi_f^\mathcal{H}(x;\theta):=
\frac{\partial_\theta
H_\theta'(x)}{H_\theta'(x)}+\frac{\partial_\theta
H_\theta(x)f'(H_\theta(x))}{f(H_\theta(x))}
\end{equation}
over $D$ (it is set to 0 outside $D$). Extracting e.c.'s from
equation (\ref{scoreeq}) is all but evident here, as (i) there is no
structural reason for $\partial_\theta H_\theta(x)$ to cross the
$x$-axis so as to allow to fix $d$ to 1 as in the
scale case over $\R$, and (ii) the generality of the model hampers a
clear understanding of the role of $\theta$ inside the
densities. Especially the latter point is crucial, as e.c.'s cannot
depend on the parameter. We thus need to further
specify the form of $f_\mathcal{H}(x;\theta)$ or, more exactly, the
form of the transformations in $\mathcal{H}$ and thus the action of
$\theta$ inside the densities. We choose to restrict our attention to
transformations $H_\theta$ satisfying the following two
factorizations:
\[
\cases{ \partial_\theta H_\theta(x)=T(\theta)U_1
\bigl(H_\theta(x)\bigr),
\vspace*{2pt}\cr
\displaystyle \frac{\partial_\theta H_\theta'(x)}{H_\theta'(x)}=T(\theta)U_2
\bigl(H_\theta(x)\bigr),}
\]
where $T$, $U_1$ and $U_2$ are real-valued functions. At first sight,
such restrictions might seem severe, but there exist numerous
one-parameter transformations enjoying these factorizations, including:
\begin{itemize}
\item[-] transformations of the form $H_\theta(x)=a_1(x)+a_2(\theta
)$ defined over the entire real line, with $a_1$ a monotone increasing
differentiable function over $\R$ and $a_2$ any real-valued
differentiable function. These transformations lead to ``generalized
location families'' and satisfy the above factorizations with $T(\theta
)=a_2'(\theta)$, $U_1(x)=1$ and $U_2(x)=0$;
\item[-] transformations of the form $H_\theta(x)=a_1(x)a_2(\theta)$
defined over $\R, \R_0^+$ or $\R_0^-$, with $a_1$ a monotone
increasing differentiable function over the corresponding domain and
$a_2$ a positive real-valued differentiable function. These
transformations lead to ``generalized scale families'' and satisfy the
above factorizations with $T(\theta)=a_2'(\theta)/a_2(\theta)$,
$U_1(x)=x$ and $U_2(x)=1$;
\item[-] transformations of the form
$H_\theta(x)=\operatorname{sinh}(\operatorname{arcsinh}(x)+\theta)$
defined over $\R$. These are the so-called sinh--arcsinh transformations
put to use in Jones and Pewsey \cite{JP09} in order to define
sinh--arcsinh distributions which allow to cope for both skewness and
kurtosis. The above factorizations are verified for $T(\theta)=1$,
$U_1(x)=\sqrt{1+x^2}$ and $U_2(x)=x/\sqrt{1+x^2}$.
\end{itemize}
Under these premisses, equation (\ref{scoreeq}) becomes
\begin{eqnarray*}
&& d \biggl(U_2\bigl(H_\theta(x)\bigr)+U_1
\bigl(H_\theta(x)\bigr)\frac{f'(H_\theta
(x))}{f(H_\theta(x))} \biggr)
\\
&&\quad = U_2\bigl(H_\theta(x)\bigr)+U_1
\bigl(H_\theta(x)\bigr)\frac{g'(H_\theta(x))}{g(H_\theta
(x))} \qquad\forall x\in D,
\end{eqnarray*}
which can be rewritten as
\[
d \biggl(U_2(x)+U_1(x)\frac{f'(x)}{f(x)}
\biggr)=U_2(x)+U_1(x)\frac
{g'(x)}{g(x)} \qquad\forall x\in S.
\]
This first-order differentiable equation admits as solution
$g(x)=c \exp((d-1)\int^xU_2(y)/\allowbreak U_1(y)\mrmd y)(f(x))^{d}$ for some $d>0$
(such that $g$ is integrable) and $c$ a normalizing constant. This
relationship establishes the $\mathcal{H}$-based e.c.'s. As for the
scale case, $c=d=1$ when there exists $x_0\in S$ such that
$U_1(x_0)=0$, yielding e.c.'s constituted of singletons $\{f\}$.

For each transformation group $\mathcal{H}$, we obtain the following MLE
characterization theorem for one-parameter group families. The proof of
this result contains nothing new and is thus omitted.

\begin{theorem}\label{TheoH}
Let $\mathcal{F}(\mathcal{H})$ and $\mathcal{G}(\mathcal{H})$ be two
distinct $\mathcal{H}$-based e.c.'s and let their respective
representatives $f$ and $g$ be two continuously differentiable
densities with common full support $S$. Let
$\varphi_{f}^{\mathcal{H}}:=U_2(x)+U_1(x)f'(x)/f(x)$ be the
$\mathcal{H}$-score function of $f$. If $\varphi_{f}^{\mathcal{H}}$ is
invertible over
$S$ then there exists
$N \in\N$ such that, for any $n\ge N$, we have
$\hat{\theta}^{(n)}_f=\hat{\theta}^{(n)}_g$ for all samples of size
$n$ if
and only if there exist constants $c,d\in\R_0^+$ such that $g(x)=c
\exp((d-1)\int^xU_2(y)/U_1(y)\mrmd y)(f(x))^d$ for all $x\in S$ (with $c=d=1$
if there exists $x_0\in S$ such that $U_1(x_0)=0$), that is, if and
only if $\mathcal{F}(\mathcal{H})=\mathcal{G}(\mathcal{H})$.
The smallest integer for which this holds is
$\mbox{MNSS}= \max\{N_f, 3\}$,
with $N_f$ the MCSS as defined in Lemma \ref{lemmalem2}.
\end{theorem}

Aside from location- and scale-based characterizations (or
variations thereof) which are already
available from Theorems \ref{Theoloc} and \ref{Theosca}, Theorem
\ref{TheoH}
allows, inter alia, to characterize asymmetric distributions
(namely the sinh--arcsinh distributions of Jones and Pewsey \cite{JP09})
with respect to
their skewness parameter.

\section{Examples}\label{Appl}

In this section, we analyze and discuss several examples of absolutely
continuous distributions in light of the findings of the previous
sections. We indicate, in each case, the corresponding MNSS. As we
shall see, we hereby retrieve a wide variety of
existing results, and obtain several new ones. We stress that, in each
case discussed below, the minimal sample size provided is optimal in
the sense that counter-examples can be constructed if the results only
are supposed to hold true for smaller sample sizes. Moreover, we
attract the reader's attention to the fact that, in some examples of
scale-based characterizations, we need the scale-identification
condition (\ref{eqsc-cond}), whereas in others it is superfluous, as
explained in Section \ref{secscale}.

For the sake of clarity, we will adopt in this section the
commonly used notations $\hat\mu^{(n)}_f$ and $\hat\sigma^{(n)}_f$ for
location and scale ML estimators.

\subsection{The Gaussian distribution}
Consider the Gaussian distribution whose MLE characterizations for both
the location and the scale parameter have been extensively discussed
in the literature. For $\phi$ the standard Gaussian density, we get
$\varphi_\phi(x)=x$ which is invertible over $\R$ and has image
$\operatorname{Im}(\varphi_\phi)=\R$. As already mentioned several times, the MLE
$\hat\mu^{(n)}_\phi$ is given by the sample arithmetic mean
$\bar{x}$. Thus, Theorems \ref{Theoloc} and \ref{Theoloc2} apply, with
$\mbox{MNSS}=3$ since $P^+_{\phi}=P^-_{\phi}=\infty$. The first corresponds
to Azzalini and Genton \cite{AG07}, Theorem 1, the second to
Teicher \cite{T61}, Theorem 1. Regarding the scale characterization for
$\hat\sigma^{(n)}_\phi=(n^{-1}\sum_{i=1}^nx_i^2)^{1/2}$, direct
calculations reveal that $\psi_\phi(x)=1-x^2$ which is
invertible over both $\R^+_0$ and $\R^-_0$ and maps both domains onto
$(-\infty,1)$. The conditions of Theorem \ref{Theosca2} are thus fulfilled
and yield that the MNSS equals $\infty$. Hence we retrieve
Teicher \cite{T61}, Theorem 3.

\subsection{The gamma distribution}
Consider the gamma distribution with tail parameter $\alpha>0$, whose
density is given by
\[
f(x)=\frac{1}{\Gamma(\alpha)}x^{\alpha-1}\exp(-x)\mathbb{I}_{(0,\infty)}(x),
\]
where $\mathbb{I}_A$ represents the indicator function of the set
$A$. The exponential density is a special case of
gamma densities obtained by setting $\alpha=1$. Gamma distributions
are not natural location families; this can be
seen, for instance, by considering the exponential case, where
$\varphi_f(x)=1$ and hence the location likelihood equations make no
sense. Within the framework of the current paper, we anyway do not
provide a
location-based MLE characterization of gamma densities, since their
support is only $\R_0^+$ instead of $\R$. On the contrary, gamma
densities allow for agreeable scale
characterizations. Indeed easy computations show that
$\varphi_f(x)=(-\alpha+1)/{x}+1$, $\psi_f(x)=\alpha-x$, which is thus
invertible over $\R_0^+$, $\operatorname{Im}(\psi_f)=(-\infty,\alpha)$ and
$\hat{\sigma}^{(n)}_f=\alpha^{-1}\bar{x}$. We can therefore use
Theorem \ref{Theosca} in combination with the scale-identification
condition (\ref{eqsc-cond}) to obtain that the gamma distribution with
shape $\alpha$ is characterizable \mbox{w.r.t.} its scale MLE $\alpha
^{-1}\bar{x}$ for an infinite MNSS. We
hereby recover Teicher \cite{T61}, Theorem 2, and the univariate case of
Marshall and Olkin \cite{MO93}, Theorem 5.1.

\subsection{The generalized Gaussian distribution}

Consider the one-parameter generalization of the normal distribution
proposed in Ferguson \cite{F62}, with density
\[
f(x)=\frac{|\gamma|\alpha^\alpha}{\Gamma(\alpha)}\exp\bigl(\alpha
\gamma x-\alpha\exp(\gamma x)\bigr),
\]
where $\alpha>0$ and $\gamma$, the additional
parameter, differs from zero (Ferguson has proved that, for
$\gamma\rightarrow0$, this density converges to the Gaussian). This
probability law is in fact strongly related to the gamma distribution,
as it is defined as $\gamma^{-1}\log(X)$ with
$X\sim\operatorname{Gamma}(\alpha)$. Now, direct calculations yield
$\varphi_f(x)=-\alpha\gamma(1-\exp(\gamma x))$, invertible over $\R$,
$\operatorname{Im}(\varphi_f)=\operatorname{sign}(\gamma)(-\alpha|\gamma|,\infty)$ and
$\hat{\mu}^{(n)}_f=\gamma^{-1}\log(n^{-1}\sum_{i=1}^n\exp(\gamma x_i))$.
Hence, from Theorem \ref{Theoloc}, we deduce that these distributions
can be characterized in terms of their location parameter, with MNSS
equal to $\infty$; we retrieve Ferguson \cite{F62}, Theorem 5.
Concerning the
scale part,
$\psi_f(x)=\alpha\gamma x(1-\exp(\gamma x))+1$ is not invertible over
the whole real line, but invertible over both $\R^+_0$ and $\R^-_0$,
and maps both half-lines onto $(-\infty,1)$. Consequently,
Theorem \ref{Theosca2} reveals that this distribution admits as well a
scale MLE characterization result, with infinite MNSS.

\subsection{The Laplace distribution}

Consider the Laplace distribution with density
\[
f(x)=\exp\bigl(-|x|\bigr)/2.
\]
One easily obtains
$\varphi_f(x)=\operatorname{sign}(x)$
and
$\psi_f(x)=-x\operatorname{sign}(x)+1$. While
the former function is clearly not invertible at all (but allows for a
location MLE characterization; see the end of
Section \ref{secscale}), the latter is invertible on both $\R_0^-$
and $\R_0^+$ with $\operatorname{Im}(\psi_f)=(-\infty,1)$. Hence
Theorem \ref{Theosca2} applies and reveals that the Laplace
distribution is also MLE-characterizable \mbox{w.r.t.} its scale
parameter (with infinite MNSS), which complements the existing results
on MLE characterizations of the Laplace distribution from Ghosh and Rao
\cite{GR71},
Kagan, Linnik and Rao \cite{KLR73}, Marshall and Olkin \cite{MO93}.

\begin{corollary} The statistic
\[
\hat\sigma^{(n)}_f= \Biggl(n^{-1}\sum
_{i=1}^n|x_i| \Biggr)^{-1}
\]
is
the MLE of the scale parameter $\sigma$ within scale families over $\R
$ for all samples of all sample sizes if and only if the samples
are drawn from a Laplace distribution.
\end{corollary}

For the sake of readability we will, here and in the sequel, content
ourselves with such informal
statements of our characterization results; rigorous statements are
straightforward adaptations of the corresponding theorems from the
previous sections.

\subsection{The Weibull distribution}

Consider the Weibull distribution with
density
\[
f(x)=kx^{k-1}\exp\bigl(-x^k\bigr)\mathbb{I}_{(0,\infty)}(x),
\]
where $k>0$ is the shape parameter. As for
gamma distributions, we do not provide a location-based MLE
characterization for this distribution on the positive real
half-line. Regarding the scale part, we have
$\varphi_f(x)=-\frac{k-1}{x}+kx^{k-1}$, $\psi_f(x)=k(1-x^k)$, clearly
invertible over $\R_0^+$, and $\operatorname{Im}(\psi_f)=(-\infty,k)$.
Thus, all conditions
for Theorem \ref{Theosca} are satisfied, from which we derive, under
the scale-identification condition (\ref{eqsc-cond}), the
following, to the best of our knowledge new, MLE characterization of
the Weibull distribution.

\begin{corollary} Let condition (\ref{eqsc-cond}) hold. Then the statistic
\[
\hat\sigma^{(n)}_f= \Biggl(n^{-1}\sum
_{i=1}^nx_i^k
\Biggr)^{-1/k}
\]
is
the MLE of the scale parameter $\sigma$ within scale families over $\R
_0^+$ for all samples of all sample sizes if and only if the samples
are drawn from a Weibull distribution with shape parameter $k$.
\end{corollary}

\subsection{The Gumbel distribution}

Consider the Gumbel distribution with
density
\[
f(x)=\exp\bigl(-x-\exp(-x)\bigr).
\]
Straightforward manipulations
yield $\varphi_f(x)=1-\exp(-x)$, invertible over $\R$ and $\operatorname{Im}(\varphi_f)=(-\infty, 1)$. Thus, all conditions
for Theorem \ref{Theoloc} are satisfied, from which we derive the
following, to the best of our knowledge new, MLE characterization of
the Gumbel distribution (actually, of the \textit{power-Gumbel distribution}).

\begin{corollary}
The statistic
\[
\hat\mu^{(n)}_f=\log\Biggl[ \Biggl( n^{-1}\sum
_{i=1}^n\exp(-x_i)
\Biggr)^{-1} \Biggr]
\]
is the MLE
of the location parameter $\mu$ within location families over $\R$
for all samples of all sample sizes if and only if the samples are drawn
from a power-Gumbel distribution with density $c\exp(-dx-d\exp(-x))$
for $c,d\in\R_0^+$.
\end{corollary}

As for the scale part, it follows that $\psi_f(x)=x(-1+\exp(-x))+1$,
non-invertible over $\R$ but invertible over both $\R_0^+$ and $\R_0^-$,
and $\operatorname{Im}(\psi_f)=(-\infty,1)$. Consequently,
Theorem \ref{Theosca2} applies and shows that the Gumbel distribution
(here not a general power-Gumbel distribution)
allows as well for a MLE characterization with respect to its scale
parameter (with corresponding MNSS equal to $\infty$).



\subsection{The Student distribution}
Consider the Student distribution with $\nu>0$ degrees of freedom,
with density
\[
f(x) = \kappa(\nu) \biggl(1+\frac{x^2}{\nu} \biggr)^{-(\nu+1)/2}
\]
for $\kappa(\nu)$ the appropriate normalizing constant. Then although the
support of $f$ is the whole real line, the location score function
$\varphi_f(x)=(\nu+1)\frac{x}{\nu+x^2}$ is\vspace*{1pt}
not invertible and thus we cannot provide a location
characterization. On the other hand straightforward computations yield
\[
\psi_f(x)=1-(\nu+1)\frac{x^2}{\nu+x^{2}}=-\nu+\frac{\nu(\nu
+1)}{\nu+x^2}.
\]
This function is invertible over both the positive and the negative
real half-line with $\operatorname{Im}(\psi_f)=(-\nu,1)$ and thus the Student
distribution
with $\nu$ degrees of freedom is by virtue of Theorem \ref{Theosca2}
scale-characterizable with
\[
\mathrm{MNSS} =\cases{ \displaystyle \biggl\lceil1+\frac{1}{\nu}\biggr\rceil, &\quad if $\nu<1$,
\vspace*{2pt}\cr
3, &\quad if $\nu=1$,
\vspace*{2pt}\cr
\lceil1+\nu\rceil, &\quad if $\nu>1$.}
\]
This result generalizes the scale characterization of the Gaussian
distribution, which is a particular case of the Student distributions
when $\nu$ tends to infinity. Note that the expression above then
indeed yields an infinite MNSS. Moreover, the Cauchy distribution,
obtained for $\nu=1$, is MLE-characterizable \mbox{w.r.t.} its scale
parameter with an MNSS of 3.

\subsection{The logistic distribution}
Consider the logistic distribution, whose density is given by
\[
f(x)=\frac{\RMe^{-x}}{(1+\RMe^{-x})^2}.
\]
Straightforward manipulations yield $\varphi_f(x)=\tanh(x/2)$, which
is invertible over $\R$ and $\operatorname{Im}(\varphi_f)=(-1,1)$.
Theorem \ref{Theoloc} applies and yields the (to the best of our
knowledge) first MLE characterization of the (power-)logistic
distribution with respect to its location parameter (with corresponding
MNSS equal to $3$).

Further, $\psi_f(x)=1-x\tanh(x/2)$ is not invertible over the whole
real line, but is invertible over both $\R_0^+$ and $\R_0^-$, and
maps both half-lines onto $(-\infty,1)$. Consequently, Theorem \ref
{Theosca2} applies and yields the (to the best of our knowledge) first
scale MLE characterization of the logistic distribution, with infinite MNSS.

\subsection{The sinh--arcsinh skew-normal distribution}

As a final example, we consider the sinh--arcsinh skew-normal
distribution of Jones and Pewsey \cite{JP09} whose density is given by
\[
f(x)=\frac{1}{\sqrt{2\uppi }}\frac{(1+\operatorname{sinh}^2
(\operatorname{arcsinh}(x)+\delta))^{1/2}}{(1+x^2)^{1/2}}\RMe^{-\operatorname{sinh}^2
(\operatorname{arcsinh}(x)+\delta)/2},
\]
where $\delta\in\R$ is a skewness parameter regulating the asymmetric
nature of the distribution. Clearly, for $\delta=0$, corresponding to
the symmetric situation, one retrieves the standard normal
distribution. Now, straightforward but tedious calculations provide us
with expressions for $\varphi_f$ and $\psi_f$ which can both be seen to
be non-invertible. Hence, no location-based nor scale-based
characterizations can be obtained.
However, the sinh--arcsinh skew-normal
distribution can be characterized $\mbox{w.r.t.}$ its skewness
parameter. As shown in Section \ref{secother-parameters}, the
sinh--arcsinh transform belongs to the class of transforms leading to
group families. Consequently, its skewness score function is given by
\[
\varphi_f^{\mathcal{H}}(x)=U_2(x)+U_1(x)
\frac{\phi'(x)}{\phi
(x)}=\frac{-x^3}{(1+x^2)^{1/2}}
\]
with $\phi$ the standard Gaussian density. This mapping is invertible
over $\R$ with symmetric image~$\R$. Theorem \ref{TheoH} therefore
applies and yields the (to the best of our knowledge) first MLE
characterization of the sinh--arcsinh skew-normal distribution (with
respect to its skewness parameter) with an MNSS equal to 3.

\section{Discussion and open problems}\label{secdiscu}

In this article, we have provided a unified treatment of the topic of
MLE characterizations for one-parameter group families of absolutely
continuous distributions
satisfying certain
regularity conditions. A natural question of
interest is then in how far our methodology can be adapted to other
distributions which do not satisfy these assumptions. Of particular
interest are (i) parametric families whose score function is
either not invertible or not differentiable at a countable number of
points (such as, e.g., the Laplace distribution w.r.t. its location
parameter), (ii) families
depending on more than one parameter and (iii) discrete families.
Although we will not cover these questions in full here, we conclude
the paper by providing a number of intuitions on these questions; in
all cases it seems clear that our methodology provides -- at least in
principle -- the path towards a satisfactory answer.

Regarding the first point, an interesting issue to investigate is how
the non-invertibility of the score function influences the MNSS.
Indeed in the case of a Laplace target the MNSS is
known to be equal to 4 (see Ghosh and Rao \cite{GR71}, Kagan, Linnik
and Rao~\cite{KLR73}). This increase
is due to the fact that the
Laplace score function only takes on two distinct non-zero values so that
having (\ref{MLEeqloc}) for three
sample points forces one of the observations to be 0
(otherwise the equality cannot hold) and therefore
the case $n=3$ provides no more information than the case $n=2$ (and
thus $\mbox{MNSS}\ge4$). It
would of course be interesting to understand the influence of the
number of distinct values taken by a given score function on the corresponding
MNSS. One would, moreover, need to deal in this case with commensurability
issues in order for the corresponding identity (\ref{MLEeqloc}) to
hold; this would most certainly lead to interesting discussions. Aside
from these issues, however, the question of characterizability is, to the
best of our understanding, covered by our approach (see the end of
Section \ref{secloc}).

Regarding the second point, it seems straightforward (but clearly
requires some care) to extend our method to a multi-dimensional
location parameter, as is already done in Marshall and Olkin \cite
{MO93} for a Gaussian
target density. In a nutshell, it suffices to project the now
multi-dimensional location score function onto distinct-directional
unit vectors and then proceed ``as in the univariate case''. On the
contrary, dealing with a high-dimensional scale
parameter seems more difficult, as the scale parameter becomes a
matrix-valued scatter or shape parameter. One possibility could be to
try to adapt Marshall and Olkin's \cite{MO93} working scheme, who have
been able to
provide an MLE characterization for the scatter parameter of a
multinormal distribution. Along these lines a final issue that we
have not considered is that of MLE
characterizations of univariate target distributions with respect to
multivariate
parameters (such as the Gaussian in terms of its two parameters $(\mu,
\sigma)$). In support of our optimism for these multivariate setups, see
Duerinckx and Ley \cite{DL12} where our methodology was successfully
applied to the
(perhaps more complex) case of the spherical location parameter families.

Finally concerning the discrete setup, it seems clear that our approach
again yields in principle a satisfactory answer for discrete group
families, although this will require a certain amount of work.
We defer the systematic treatment of this
interesting
question to later publications.

We conclude this paper by an intriguing question suggested by an
anonymous referee, which is the following: do characterization results
survive mixtures of distributions? More concretely, if a given
distribution is MLE-characterizable \mbox{w.r.t.} a given parameter of
interest, under which conditions will mixtures of this distribution
remain MLE-characterizable? This question is all but straightforward to
answer. Indeed, we have seen in Section \ref{Appl} that the Student
distribution is only characterizable \mbox{w.r.t.} its scale
parameter, whereas the logistic distribution admits an MLE
characterization \mbox{w.r.t.} both the location and scale parameter;
the Student and the logistic distribution are both scale mixtures of
the Gaussian distribution, which is MLE-characterizable both as a
location and scale family. We leave this question as an open problem.

\section*{Acknowledgements}

All three authors thank Johan Segers and Davy Paindaveine for their
interesting comments during presentations of preliminary versions of
this paper. The authors are also very grateful to two anonymous
referees for extremely interesting remarks and suggestions that have
led to a definite improvement of the present paper.

Christophe Ley thanks the Fonds National de la Recherche Scientifique,
Communaut\'{e} fran\c{c}aise de Belgique, for support via a Mandat de
Charg\'{e} de Recherche FNRS. Christophe Ley is also a member of ECARES.


%

\printhistory

\end{document}